\documentclass{article}
\usepackage{amsmath}
\usepackage{pstricks}
\usepackage{pst-node}
\usepackage{amsthm}
\usepackage{amssymb}
\newcommand{\ints}{\mathbb{Z}}
\newcommand{\vect}[1]{\mathbf{#1}}
\newcommand{\EY}{{\overline Y}}\newcommand{\EU}{{\overline U}}
\newcommand{\EW}{{\overline W}}\newcommand{\Echi}{{\overline\chi}}
\newcommand{\Ep}{{\overline p}}\newcommand{\Et}{{\overline\tau}}
\newcommand{\Em}{{\overline m}}\newcommand{\Ex}{{\overline{\vect x}}}
\newtheorem{proposition}[equation]{Proposition}
\newtheorem{corollary}[equation]{Corollary}
\newtheorem{theorem}[equation]{Theorem}

\newlength\PullBackLength
\newcommand\PullBack[1][2]{%
  \setlength{\global\PullBackLength}{#1em}%
  \kern\PullBackLength%
  &
  \kern-\PullBackLength}
\makeatletter
\@addtoreset{equation}{section}
\makeatother

\begin{document}

\title{The Equivalence of Two Graph Polynomials And a Symmetric Function}%
\author{Criel Merino\thanks{Supported by CONACYT of M\'exico.}\\
Instituto de Matem\'aticas\\ Universidad Nacional Aut\'onoma de M\'exico\\
      Area de la investigaci\'on cient\'{\i}fica\\ Circuito Exterior, C.U. Coyoac\'an 04510,\\
M\'exico, D.F. M\'exico.\\[5mm] Steven D. Noble
\thanks{Partially supported by the Heilbronn Institute for
Mathematical Research, Bristol, U.K. and by the Royal Society.}\\Department of Mathematical Sciences\\Brunel University\\
Kingston Lane\\Uxbridge\\UB8 3PH\\United Kingdom}%
\maketitle

\begin{abstract}
The $U$-polynomial, the polychromate and the symmetric function
generalization of the Tutte polynomial due to Stanley are known to
be equivalent in the sense that the coefficients of any one of
them can be obtained as a function of the coefficients of any
other. The definition of each of these functions suggests a
natural way in which to generalize them which also captures
Tutte's universal $V$-functions as a specialization. We show that
the equivalence remains true for the extended functions thus
answering a question raised by Dominic Welsh.
\end{abstract}


\maketitle


\section{Introduction}
This paper answers a question posed by Dominic Welsh in a talk in
2005~\cite{welsh:talk05} concerning the notions of equivalence and
specialization of graph polynomials and symmetric functions. We
say that for graph polynomials $P$ and $Q$, $P$ \emph{specializes}
to $Q$ written $P \succ Q$ if the coefficients of $Q$ may be
obtained as functions of the coefficients of $P$ and the number of
vertices of the graph. Graph polynomials $P$ and $Q$ are
\emph{equivalent} if $P \succ Q$ and $Q\succ P$. These notions may
be extended to symmetric functions by allowing the coefficients to
be those of a symmetric function with respect to some basis.
Defining equivalence in the right way is not completely
straightforward~\cite{mak:zoo}, but this very simplistic notion
will suffice for our purposes. Equivalence of graph polynomials or
symmetric functions is clearly an equivalence relation.

Many specializations of the Tutte polynomial are well-known and
include the chromatic and reliability polynomials. The key objects
in this paper are extensions of two graph polynomials and one
symmetric function that themselves generalize the Tutte
polynomial. Brylawski introduced the
polychromate~\cite{bry:intersection}, which is a polynomial in
countably infinitely many variables. Stanley generalized the
definition of the chromatic polynomial~\cite{stanley:symm} to the
chromatic symmetric function and a little later extended this to a
symmetric function generalization of the Tutte
polynomial~\cite{stanley:garsia}. For brevity we call this the
Tutte symmetric function. Motivated by problems from knot theory,
Noble and Welsh introduced the
$U$-polynomial~\cite{noble+welsh:U+W} and proved that it is
equivalent to the Tutte symmetric function.
Sarmiento~\cite{sarm:poly} then showed that the $U$-polynomial and
the polychromate are equivalent.

Taking for the moment an entirely naive and informal perspective,
the definitions of these three functions lack a certain symmetry.
In each of them the $x$ variable of the Tutte polynomial is
generalized to countably infinitely many variables whereas the $y$
variable remains essentially unchanged. More importantly none of
the three functions specializes to Tutte's universal
$V$-functions~\cite{tutte:ring}. It turns out that there are
natural ways to define extended versions of each of these
polynomials to overcome this problem which at the same time
address the lack of symmetry.

The question asked by Welsh~\cite{welsh:talk05} was whether the
equivalence of the $U$-polynomial and the polychromate carries
over to their extended versions. Sarmiento's proof is quite
involved and extending her methods did not appear to be an easy
task. A key step in our approach is to introduce an extension of
the Tutte symmetric function as an intermediate object between the
two polynomials. Our main results are that the extended Tutte
symmetric function is equivalent to both the extended polychromate
and the extended $U$-polynomial. Since equivalence is transitive
this answers Welsh's question. Our proof suggests a way to
simplify Sarmiento's proof.

An overview of the paper is as follows. In the next section we
present some preliminary definitions concerning symmetric
functions. We then define the previously studied polynomials
covered here and briefly survey some of their properties. Section
four contains the definitions of two new graph polynomials and one
new graph symmetric function together with our main results. We
end with a brief conclusion and an open problem.

\section{Partitions and symmetric functions}
We give some definitions and fix some notation which we will use
throughout the paper. Note that all of our graphs are finite and
may have multiple edges and loops. When the underlying graph is
obvious we use $V$ and $E$ to denote respectively its set of
vertices and edges and let $n=|V|$ and $m=|E|$. Given a graph $G$,
$G|A$ is formed by deleting all the edges in $E\setminus A$ (but
keeping all the vertices). We use $k(G)$ to denote the number of
connected components of $G$ and define the rank of a set $A$ of
edges to be given by $r(A)=|V|-k(G|A)$. If $A\subseteq E$ then let
$\pi(A)$ denote the partition of $V$ for which the blocks are the
connected components of $G|A$.

Given a partition $\pi$ of a set $A$, its \emph{type} is the
integer partition of $|A|$ for which the parts are the sizes of
the blocks of $\pi$. If $\tau$ is an integer partition of $n$, we
write $\tau \vdash n$ and let $k(\tau)$ be the number of parts of
$\tau$. As usual we write the components of an integer partition
in decreasing order so that if $\tau = (n_1,n_2,\ldots,n_k)$, we
have $n_1\geq n_2 \geq \ldots \geq n_k$.

We next introduce two symmetric function bases. For $r \geq 1$ let
\[ p_r(\vect x) = p_r (x_1,x_2,\ldots) = \sum_{i=1}^{\infty}x_i^r.\]
Now suppose that $\tau = (n_1,\ldots,n_k)$ is an integer partition
of $n$. Then we define $p_\tau(\vect x)$ to be the symmetric
function $\prod_{i=1}^k p_{n_i}(\vect x)$. The collection
$\{p_\tau(\vect x): \tau \vdash n\}$ forms a basis called the
\emph{power-sum basis} for the homogenous symmetric functions of
degree $n$ in $\vect x$~\cite{macdonald:symm}.

Elements of the second basis, the \emph{augmented monomial basis}
are also defined in terms of an integer partition
$\tau=(n_1,\ldots,n_k)$ of $n$. Let
\[ m_\tau(\vect x) = \sum_{(i_1,\ldots,i_k)} x_{i_1}^{n_1} \cdots
x_{i_k}^{n_k},\] where the sum is over all $k$-tuples of pairwise
distinct strictly positive integers. Again, the collection
$\{m_\tau(\vect x) : \tau \vdash n\}$ is a basis for the
homogenous symmetric functions of degree $n$ in $\vect x$. Note
that $m_\tau$ is often denoted by $\tilde m_\tau$.

We now generalize some of these ideas to what we call paired
symmetric functions. We have not been able to find any reference
to these objects in the literature but surely they have been
encountered many times before. First we define an \emph{integer
pair partition} of a pair of strictly positive integers $(a,b)$ to
be a list of pairs of integers $((a_1,b_1),\ldots,(a_k,b_k))$ such
that
\begin{enumerate}
\item
$(a_1,\ldots,a_k)$ is an integer partition  of $a$;
\item
for all $i$, $b_i$ is a non-negative integer and $\sum_{i=1}^k b_i
\leq b$;
\item if $i<j$ then either $a_i>a_j$ or $a_i=a_j$ and $b_i\geq b_j$, that is the pairs $(a_i,b_i)$ are written in lexicographically
decreasing order.
\end{enumerate}
If $\Et$ is an integer pair partition of $(a,b)$, we write $\Et
\vdash (a,b)$.

The canonical example of an integer pair partition is as follows.
For a graph $G$, let $\pi$ be a partition of its vertex set. Order
the blocks of $\pi$ in an arbitrary way. Let $a_i$ denote the
number of vertices in the $i$th block and let $b_i$ denote the
number of edges of $G$ having both endpoints in the $i$th block.
Now write the pairs $(a_i,b_i)$ in lexicographically decreasing
order to obtain the integer pair partition which we denote
$\Et(\pi)$.

We next define a paired symmetric function. Suppose $f$ is a
function in the pairs of variables
$(x_1,t_1),\ldots,(x_i,t_i),\ldots$ such that for any permutation
$\pi$ of $\ints^{>0}$
\[f((x_{\pi(1)},t_{\pi(1)}),\ldots,(x_{\pi(i)},t_{\pi(i)}),\ldots)
=f((x_{1},t_{1}),\ldots,(x_i,t_i),\ldots).\] We require
additionally that $f$ is homogenous in the $\vect x$ variables.
Then we call $f$ a \emph{paired symmetric function}. Notice that
$f$ is not generally a symmetric function in the usual sense. A
key observation is that it is possible to extend the two classes
of symmetric function bases discussed above to paired symmetric
functions.

We describe first how to extend the definition of the power-sum
basis. If $r\in \ints^{>0}$ and $s\in \ints^{\geq 0}$ define
$\Ep_{r,s}(\vect x,\vect t) = \sum_{i=1}^{\infty} x_i^r t_i^s$.
For an integer pair partition $\Et$ let
\[\Ep_{\Et}(\vect x,\vect t) = \prod_{(a_i,b_i)\in \Et} \Ep_{a_i,b_i}(\vect x,\vect t).\]
Then the collection $\{\Ep_{\Et}(\vect x,\vect t) : \Et \vdash
(n,m)\}$ forms a basis for the paired symmetric functions of
degrees $n$ and $m$ in respectively $\vect x$ and $\vect t$.

A second basis for the paired symmetric functions is defined by
extending the definition of the augmented monomial basis. If
$\Et=((a_1,b_1),\ldots,(a_k,b_k))$ is an integer pair partition of
$(n,m)$ then let
\[ \Em_{\Et}(\vect x,\vect t) = \sum_{i_1,\ldots,i_k} x_{i_1}^{a_1} (1+t_{i_1})^{b_1}
\cdots x_{i_k}^{a_k} (1+t_{i_k})^{b_k},\] where the summation is
over all $k$-tuples $(i_1,\ldots,i_k)$ of pairwise distinct
strictly positive integers. The collection $\{\Em_{\Et}(\vect
x,\vect t) : \Et \vdash (n,m)\}$ forms a basis for the paired
symmetric functions of degrees $n$ and $m$ in respectively $\vect
x$ and $\vect t$.

\section[title]{A menagerie of polynomials
\footnote{The section title is suggested by the
    title of~\cite{mak:zoo}}}
We give definitions of and some relations between a collection of
graph polynomials beginning with two very well-established
examples and moving on to four that are more recent.

The chromatic polynomial $P_G(\lambda)$ was introduced by Birkhoff
in 1912~\cite{birkhoff:1912} in an effort to prove the four colour
theorem. For a positive integer $\lambda$ it is defined to be the
number of proper colourings of the vertices of $G$ using colours
drawn from a set of size $\lambda$. Whitney~\cite{whitney32}
showed that
\begin{equation}\label{eq:whitney}
P_G(\lambda) = \sum_{A \subseteq E} (-1)^{|A|} \lambda^{k(G|A)}.
\end{equation}
This relation is one way to prove that the chromatic polynomial is
truly a polynomial but more importantly a generalization of it
forms the crux of one of our proofs.

Tutte introduced his eponymous polynomial in~\cite{tutte:ring}. Given a graph
$G$, the Tutte polynomial $T_G(x,y)$ is given by
\[
T_G(x,y) = \sum_{A \subseteq E} (x-1)^{r(E)-r(A)} (y-1)^{|A|-r(A)}.
\]
Using \eqref{eq:whitney} one obtains the well-known specialization
\[ P_G(\lambda) = \lambda^{k(G)} T_G(1-\lambda,0).\]
Because $k(G)=|V|-r(E)$ and $r(E)$ is determined by the
coefficients of the Tutte polynomial we obtain $T \succ P$.

The Tutte polynomial contains a whole host of specializations, for
example the number of spanning trees, number of spanning forests
and the reliability polynomial as well as applications in
statistical mechanics and knot theory. Details of many
specializations are contained in~\cite{bry+ox:tutte,welsh:lmslns}.

Motivated by a series of papers~\cite{chmutov:1,chmutov:2,chmutov:3},
the weighted graph polynomial $U$ was introduced
in~\cite{noble+welsh:U+W}. The authors
of~\cite{chmutov:1,chmutov:2,chmutov:3} introduce a graph polynomial
derived from Vassiliev invariants of knots and note that this
polynomial does not include the Tutte polynomial as a special
case. With a slight
generalisation of their definition we obtain the weighted graph
polynomial $U$ that does include the Tutte polynomial.

The original definition of $U$ involved a recurrence relation
using deletion and contraction, but for the purposes of this paper
it is most useful to define $U$ using the ``states model
expansion'' from Proposition 5.1 in~\cite{noble+welsh:U+W}.
\begin{equation}\label{eq:state}
U_G(\vect x,y)=U_G(x_1,x_2,\ldots,y) = \sum_{A \subseteq E}
x_{n_1} x_{n_2} \cdots x_{n_{k(G|A)}}
(y-1)^{|A|-r(A)},\end{equation} where $n_1,\,\ldots,n_{k(G|A)}$
are the numbers of vertices in the connected components of $G|A$
and $x_1,x_2,\ldots$ are commuting indeterminates. For example, if
$G$ is a triangle then
\[U_G(\vect x,y) = x_1^3 + 3x_1x_2 + 3x_3 + (y-1)x_3 =
x_1^3+3x_1x_2+2x_3+yx_3.\]

The next few results are all proved in~\cite{noble+welsh:U+W}.
The first result shows that $U \succ T$.
\begin{proposition}
For any graph $G$,
\[ T_G(x,y) = (x-1)^{-k(G)}U_G(x_i=x,y).\]
\end{proposition}
Note that we have abused notation somewhat by writing
$U_G(x_i=x,y)$ where we mean setting $x_i=x$ for all $i$.

The attraction of $U$ is that it contains many other graph
invariants as specialisations, for instance the $2$-polymatroid
rank generating function of Oxley and
Whittle~\cite{oxley+whittle:2-polymatroid}, and as a consequence
the matching polynomial and the stable set polynomial~\cite{farr}.

A \emph{stable} set in a graph $G$ is a set $S$ of vertices for which
$G$ has no edge with both endpoints in $S$. The stability polynomial
$A_G(p)$ was introduced by Farr in~\cite{farr} and is given by
\[ A_G(p) = \sum_{U \in \mathcal{S}(G)} p^{|U|}(1-p)^{|V(G)\setminus
  U|},\]
where $\mathcal {S}(G)$ is the set of all stable sets of $G$.
\begin{proposition}
If $G$ is loopless then $A(G;p)$ is given by
\[ A_G(p)=U_G(x_1=1,\,x_j=-(-p)^j \textrm{ for }j \geq 2 ,\,y=0).\]
\end{proposition}

The two-polymatroid rank generating function $S_G(u,v)$ was
introduced by Oxley and Whittle
in~\cite{oxley+whittle:2-polymatroid} and is defined as follows.
Given a graph $G$ and $A \subseteq E(G)$ let $f(A)$ denote the
number of vertices of $G$ that are an endpoint of an edge in $A$.
Then
\[ S_G(u,v) = \sum_{A \subseteq E(G)} u^{|V(G)|-f(A)}v^{2|A|-f(A)}.\]
$S$ contains the matching polynomial as a specialisation.
\begin{proposition}
Let $G$ be a loopless graph with no isolated vertices. Then
\[ S_G(u,v) = U_G(x_1=u,\,x_2=1,\,x_j=v^{j-2} \textrm{ for }j >
2,\,y=v^2+1).\]
\end{proposition}

The chromatic symmetric function was
developed by Stanley in~\cite{stanley:symm}. Let $G$ be a graph with
vertex set $V=\{v_1,\ldots,v_n\}$. Then $X_G$ is a homogeneous
symmetric function of degree $n$ defined by
\[X_G(\vect x) = X_G(x_1,x_2,\ldots) = \sum_{\chi} x_{\chi(v_1)}  x_{\chi(v_2)}
\cdots x_{\chi(v_n)},\] where the sum ranges over all proper
colourings $\chi:V \rightarrow \ints^{>0}$.

The following result from~\cite{noble+welsh:U+W} shows that
$U\succ X$.
\begin{proposition}
For any graph $G$
\[ X_G(\vect x)= (-1)^{|V|}U_G(x_j=-p_j,y=0).\]
\end{proposition}

In a second paper on the chromatic symmetric
function~\cite{stanley:garsia}, Stanley introduced the Tutte
symmetric function defined by
\[Y_G(\vect x,t) = Y_G(x_1,x_2,\ldots,t) = \sum_{\chi} x_{\chi(v_1)}  x_{\chi(v_2)}
\cdots x_{\chi(v_n)} (1+t)^{b(\chi)},\] where the sum is now over
all colourings $\chi:V \rightarrow \ints^{>0}$ and $b(\chi)$ is
the number of monochromatic edges, that is, edges for which both
endpoints receive the same colour.

In~\cite{noble+welsh:U+W}, the following was shown.
\begin{theorem}\label{th:Usymm}
The polynomial $U$ and symmetric function $Y$ are equivalent. In
particular $Y_G $ is easily obtained from $U_G$ by the
substitution
\[Y_G(\vect x,t)= t^{|V|}U_G\left(x_j=\frac{p_j(\vect x)}t,y=t+1\right).\]
Conversely, if we expand $Y_G$ in terms of the power-sum basis
then we can recover $U_G$.
\end{theorem}
If $\tau=(n_1,\ldots,n_k)$ then we use $\vect x_\tau$ to denote
the monomial $\prod_{i=1}^k x_i^{n_i}$. Another way of describing
the substitution into $U_G$ that produces $Y_G$ is to say that if
for each $\tau \vdash n$ and each $i$ the monomial $\vect x_\tau
y^i$ in $U_G$ is replaced by $p_\tau(\vect x)
t^{n-k(\tau)}(t+1)^i$ then $Y_G$ is obtained.

The final polynomial that we will define is the polychromate,
introduced originally by Brylawski~\cite{bry:intersection}. Given
a graph $G$ and a partition $\pi$ of its vertices into non-empty
blocks, we define $e(\pi)$ to be the number of edges with both
endpoints in the same block of the partition.

The polychromate $\chi_G(\vect x,y)$ is defined by
\[ \chi_G(\vect x,y) = \sum_{\pi} y^{e(\pi)} \vect x_{\tau(\pi)},\]
where the summation is over all partitions of $V(G)$.

The following result is due to Sarmiento~\cite{sarm:poly}.
\begin{theorem}
The polynomials $U$ and $\chi_G$ are equivalent.
\end{theorem}
Obtaining $U$ from $\chi_G$ or vice versa is complicated and we do
not explain this here but discuss it further at the end of the
next section.

\section{Extensions}
\subsection{The extended Tutte symmetric function}
Our extension of Stanley's Tutte symmetric function replaces the
$t$ variable by countably infinitely many variables
$t_1,t_2,\ldots,$ enumerating not just the total number of
monochromatic edges but the numbers of monochromatic edges of each
colour. It is defined as follows. \begin{equation}\EY_G(\vect
x,\vect t) = \sum_{\chi}\left(\prod_{i=1}^n x_{\chi(v_i)}\right)
\left(\prod_{i=1}^\infty (1+t_i)
^{b_i(\chi)}\right),\label{eq:extsymmdef}\end{equation} where the
sum is over all colourings $\chi:V \rightarrow \ints^{>0}$ and
$b_i(\chi)$ is the number of monochromatic edges for which both
endpoints have colour $i$.

The function $\EY$ is a paired symmetric function of degrees $n$
and $m$ in respectively the $\vect x$ and $\vect t$ variables.
Note that $\EY_G$ is not homogenous in $\vect t$ unless $m=0$.

We can obtain a version of \eqref{eq:whitney} which applies to the
extended Tutte symmetric function.
\begin{proposition}\label{prop:mobius}
For any graph $G$,
\[\EY_G(\vect x, \vect t) = \sum_{A\subseteq E} \Ep_{\Et(\pi(A))}(\vect x,\vect t).\]
\end{proposition}
\begin{proof}
Given a colouring $\chi$, let $B_i(\chi)$ denote the monochromatic
edges for which both endpoints have colour $i$. Furthermore let
$B(\chi)=\bigcup_i B_i(\chi)$, the set of all monochromatic edges.
For each $i$ we can write
\[ (1+t_i)^{b_i(\chi)} = \sum_{A_i\subseteq B_i(\chi)} t_i^{|A_i|}.\]
So we have
\begin{align*}
\EY_G(\vect x,\vect t) &= \sum_{\chi}\left(\prod_{i=1}^n
x_{\chi(v_i)}\right) \left(\prod_{i=1}^\infty \sum_{A_i\subseteq
B_i(\chi)} t_i^{|A_i|}\right)\\
&= \sum_{\chi}\left(\prod_{i=1}^n x_{\chi(v_i)}\right)
\sum_{A\subseteq B(\chi)} \left(\prod_{i=1}^\infty t_i^{|A\cap
B_i(\chi)|}\right).
\end{align*}
By interchanging the order of summation, we obtain
\[\EY_G(\vect x,\vect t)= \sum_{A\subseteq E} \sum_{\chi:B(\chi) \supseteq
A}\left(\prod_{i=1}^n x_{\chi(v_i)}\right)
\left(\prod_{i=1}^\infty t_i^{|A\cap B_i(\chi)|}\right).\] The
colourings appearing in the inner summation are precisely those
which are monochromatic on the edges of $G|A$. So in any such
colouring the vertices of a component of $G|A$ must all receive
the same colour and the colours of the monochromatic edges counted
in the final product are determined by the colour of the component
of $G|A$ to which they belong. Hence for any $A\subseteq E$
\[\sum_{\chi:B(\chi) \supseteq
A}\left(\prod_{i=1}^n x_{\chi(v_i)}\right)
\left(\prod_{i=1}^\infty t_i^{|A\cap B_i(\chi)|}\right) =
\Ep_{\Et(\pi(A))}(\vect x,\vect t)\] and the result follows.
\end{proof}

\subsection{The extended $\mathbf{U}$-polynomial} The extended $U$
polynomial, which we denote by $\EU$, is a polynomial in countably
many commuting variables $z_{i,j}$ where $i\in \ints^{>0}$ and
$j\in \ints^{\geq 0}$. The definition is a natural generalization
of \eqref{eq:state} and was first introduced
in~\cite{welsh:talk05}.
\begin{equation}\label{eq:extstate}
\EU_G(\vect{z})= \sum_{A\subseteq E}
z_{c_1,e_1-c_1+1}z_{c_2,e_2-c_2+1}\cdots
z_{c_{k(G|A)},e_{k(G|A)}-c_{k(G|A)}+1},
\end{equation}
where $c_i$ and $e_i$ are respectively the number of vertices and
edges in the $i$th connected component of $G|A$.

For example if $G$ is a triangle then
\[\EU_G(\vect{z})= (z_{1,0})^3 +
3z_{1,0}z_{2,0} + 3z_{3,0} + z_{3,1}.\] Observe that
\[ U_G(\vect x,y) = \EU_G(z_{ij} = x_i(y-1)^j)\]
and so $\EU \succ U$. If we take $G_1 (G_2)$ to be a path of
length two with a loop attached at a vertex of degree one (two)
then \[U_{G_1}(\vect x,y) = U_{G_2}(\vect x,y) =
y(x_3+2x_2x_1+x_1^3).\] However
\[ \EU_{G_1}(\vect z)= z_{3,1} + z_{3,0} + z_{2,1}z_{1,0} + z_{2,0}z_{1,1} + 2z_{2,0}z_{1,0} + z_{1,0}^2z_{1,1} +
z_{1,0}^3\] but
\[\EU_{G_2}(\vect z)= z_{3,1} + z_{3,0} +
2z_{2,1}z_{1,0} + 2z_{2,0}z_{1,0} + z_{1,0}^2z_{1,1} +
z_{1,0}^3.\] Unfortunately we do not know of a pair of loopless
graphs for which $U_{G_1}= U_{G_2}$ but $\EU_{G_1}\ne \EU_{G_2}$.

We now show that $\EU$ and $\EY$ are equivalent.
\begin{theorem}
The polynomial $\EU$ and the extended Tutte symmetric function are
equivalent. More precisely
\[ \EY_G(\vect x,\vect t) = \EU_G(z_{i,j} = \Ep_{i,i+j-1}(\vect
x,\vect t)).\] Furthermore if we express $\EY_G$ in terms of the
power-sum basis as
\[ \EY_G(\vect x,\vect t) = \sum_{\Et\vdash (n,m)} a_\Et \Ep_\Et(\vect x,
\vect t),\] we obtain $\EU_G$ by replacing $\Ep_{r_1,s_1} \cdots
\Ep_{r_k,s_k}$ by $z_{r_1,s_1-r_1+1} \cdots z_{r_k,s_k-r_k+1}$.
\end{theorem}
\begin{proof}
The result follows easily from Proposition~\ref{prop:mobius}. Note
that
\[ \EY_G(\vect x,\vect t) = \sum_{A\subseteq E} \Ep_{\Et(\pi(A))}(\vect
x,\vect t) = \sum_{A\subseteq E} \Ep_{c_1,e_1} \cdots
\Ep_{c_k(G|A),e_k(G|A)},\] where $c_i,e_i$ are respectively the
number of vertices and number of edges in the $i$th component of
$G|A$. Comparing this expression with~\eqref{eq:extstate} gives
the result.
\end{proof}

There is a recurrence relation for the extended $U$-polynomial
involving deletion and contraction just as there is for $U$
itself~\cite{noble+welsh:U+W}. The recurrence relation for $U$
involves a more general polynomial $W$ defined on graphs where the
vertices have strictly positive integer weights. In order to
describe the recurrence for the extended $U$-polynomial we need to
define an extended version of $W$. We use the notation
$(G,\omega)$ to describe a graph $G$ with a strictly positive
integer weight $\omega(v)$ attached at each vertex $v$. We then
let $\EW_{(G,\omega)}(\vect z)$ be a polynomial in countably many
commuting variables $z_{i,j}$ where $i\in \ints^{>0}$ and $j\in
\ints^{\geq 0}$ and be given by
\begin{equation}\label{eq:extW}
\EW_{(G,\omega)}(\vect z)= \sum_{A\subseteq E}
z_{w_1,e_1-c_1+1}z_{w_2,e_2-c_2+1}\cdots
z_{w_{k(G|A)},e_{k(G|A)}-c_{k(G|A)}+1},
\end{equation}
where $c_i$, $e_i$ and $w_i$ are respectively the number of
vertices, the number of edges and the sum of the weights on the
vertices in the $i$th connected component of $G|A$.

For example if $(G,\omega)$ is a triangle for which the vertices
have weights $a$, $b$ and $c$ then
\[\EW_{(G,\omega)}(\vect{z})= z_{a,0}z_{b,0}z_{c,0} +
z_{a,0}z_{b+c,0} + z_{b,0}z_{c+a,0} +z_{c,0}z_{a+b,0}
+3z_{a+b+c,0} + z_{a+b+c,1}.\]

If we set $\omega(v)=1$ for all $v$ then in~\eqref{eq:extW} we
have $w_i=c_i$ for each $i$ and we obtain $\EU_G(\vect z)
=\EW_{(G,\omega)}(\vect z)$.

We now define deletion and contraction of edges in a weighted
graph. For any edge $e$ of a weighted graph $(G,\omega)$, the
\emph{deletion} of $e$, denoted by $(G,\omega)-e$ is formed by
removing $e$ from $E(G)$. For a non-loop edge $e$ with endpoints
$u$ and $v$, the \emph{contraction} of $e$, denoted by
$(G,\omega)/e$ is formed by removing $e$ from $E(G)$ and
identifying the vertices $u$ and $v$ to form a new vertex $w$
having weight $\omega(u)+\omega(v)$. So both operations conserve
the total weight of the vertices.

\begin{theorem}\mbox{ }
\begin{enumerate}
\item Suppose that the only edges of $(G,\omega)$ are loops. Let
$V=\{v_1,\ldots,v_n\}$ and suppose that $\omega(v_i)=w_i$ and that
there are $e_i$ loops attached at $v_i$. Then \begin{equation}
\EW_{(G,\omega)} = \prod_{i=1}^n \sum_{j=0}^{e_i} \binom{e_i}{j}
z_{w_i,j}.\label{eq:W1}\end{equation}
\item If $e$ is an edge of $(G,\omega)$ that is not a loop then
\begin{equation} \EW_{(G,\omega)} = \EW_{(G,\omega)-e} +
\EW_{(G,\omega)/e}.\label{eq:W2}\end{equation}
\end{enumerate}
\end{theorem}
\begin{proof}
The first part follows immediately from the definition. To prove
the second part, suppose that $e$ is an edge of $(G,\omega)$ that
is not a loop. By splitting the sum in the definition of $\EW$
depending on whether or not $A$ contains $e$ we obtain
\begin{align}\EW_{(G,\omega)}(\vect z)= &\sum_{A\subseteq E}
z_{w_1,e_1-c_1+1}\cdots
z_{w_{k(G|A)},e_{k(G|A)}-c_{k(G|A)}+1}\nonumber\\
= &\sum_{A\subseteq E-e} z_{w_1,e_1-c_1+1}\cdots
z_{w_{k(G|A)},e_{k(G|A)}-c_{k(G|A)}+1}\nonumber\\
& {} + \sum_{e\in A\subseteq E} z_{w_1,e_1-c_1+1}\cdots
z_{w_{k(G|A)},e_{k(G|A)}-c_{k(G|A)}+1}\label{eq:delctrct}.
\end{align}
From now on we just write $G$ rather than $(G,\omega)$. The first
term is $\EW_{G-e}(\vect z)$ and we claim that the second term is
$\EW_{G/e}(\vect z)$. To show this we compare the terms appearing
in the second sum in~\eqref{eq:delctrct} with those in the
definition of $\EW$ applied to $G/e$. Let $A\subseteq E-e$.
Compare the connected components of $G|(A\cup e)$ and $(G/e)|A$.
One component $C$ of $G|(A\cup e)$ contains $e$. Suppose the
endpoints of $e$ are $v$ and $w$. Then there is a component of
$(G/e)|A$ for which the vertices are those of $C-\{v,w\}$ together
with the new vertex formed when $e$ was contracted. The weight of
the new vertex is $\omega(v)+\omega(w)$ and the weight of all the
other vertices in $C$ is the same in $G|(A\cup e)$ as in $(G/e)|A$
so the total weight of the component is unchanged. Since $e$ has
been removed there is one more edge in this component in $G|(A\cup
e)$ compared with $(G/e)|A$. Similarly there is one more vertex in
the is component in $G|(A\cup e)$ compared with $(G/e)|A$. Every
other component other than $C$ has the same vertices with the same
weights and the same edges in both $G|(A\cup e)$ and $(G/e)|A$.
Hence the terms appearing in the second sum in~\eqref{eq:delctrct}
are exactly those appearing in $\EW_{G/e}(\vect z)$ and so the
claim and hence the theorem are proved.
\end{proof}

To illustrate this theorem we show how to compute $\EU$ for the
following graph. \psset{labelsep=0.2}
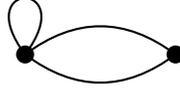
\begin{center}
\pspicture(-0.2,-0.5)(3,0.8) \cnode*(0,0){0.12}{a}
\cnode*(2,0){0.12}{b} \psbezier(0,0)(0.67,0.5)(1.33,0.5)(2,0)
\psbezier(0,0)(0.67,-0.5)(1.33,-0.5)(2,0)
\psbezier(0,0)(0.75,1)(-0.75,1)(0,0)
\endpspicture
\end{center}
To do this we add weight one to each vertex and compute $\EW$ of
the corresponding weighted graph. We use the convention that a
depiction of a graph means $\EW$ of that graph.
\begin{center}
\pspicture(-0.5,-3)(12,0.8) \cnode*(0,0){0.12}{a}
\cnode*(1.5,0){0.12}{b} \psbezier(0,0)(0.5,0.4)(1,0.4)(1.5,0)
\psbezier(0,0)(0.5,-0.4)(1,-0.4)(1.5,0)
\psbezier(0,0)(0.75,1)(-0.75,1)(0,0) \uput[d](0,0){$1$}
\uput[d](1.5,0){$1$} \rput(2.5,0){$=$}

\cnode*(3.5,0){0.12}{c} \cnode*(5,0){0.12}{d} \ncLine{-}{c}{d}
\psbezier(3.5,0)(4.25,1)(2.75,1)(3.5,0)\uput[d](3.5,0){$1$}
\uput[d](5,0){$1$} \rput(6,0){$+$} \cnode*(7.5,0){0.12}{e}
\psbezier(7.5,0)(8.5,0.75)(8.5,-0.75)(7.5,0)
\psbezier(7.5,0)(6.5,0.75)(6.5,-0.75)(7.5,0) \uput[d](7.5,0){$2$}

\rput(2.5,-1.5){$=$} \cnode*(3.5,-1.5){0.12}{f}
\cnode*(5,-1.5){0.12}{g}
\psbezier(3.5,-1.5)(4.25,-0.5)(2.75,-0.5)(3.5,-1.5)
\uput[d](3.5,-1.5){$1$} \uput[d](5,-1.5){$1$} \rput(6,-1.5){$+$}
\cnode*(7,-1.5){0.12}{h}
\psbezier(7,-1.5)(7.75,-0.5)(6.25,-0.5)(7,-1.5)
\uput[d](7,-1.5){$2$} \rput(8,-1.5){$+$}
\cnode*(9.5,-1.5){0.12}{e}
\psbezier(9.5,-1.5)(10.5,-0.75)(10.5,-2.25)(9.5,-1.5)
\psbezier(9.5,-1.5)(8.5,-0.75)(8.5,-2.25)(9.5,-1.5)
\uput[d](9.5,-1.5){$2$}

\rput(2.5,-2.5){$=$} \rput(4.25,-2.5){$z_{1,1}z_{1,0}+z^2_{1,0}$}
\rput(6,-2.5){$+$} \rput(7,-2.5){$z_{2,1}+z_{2,0}$}
\rput(8,-2.5){$+$} \rput(9.7,-2.5){$z_{2,2}+2z_{2,1}+z_{2,0}$.}
\endpspicture
\end{center}

We can now justify our claim in the introduction that the extended
$U$-polynomial (and as a corollary of the other results in this
paper the extended polychromate and extended Tutte symmetric
function) specializes to Tutte's universal $V$-functions.

The universal $V$-function is a polynomial in the commuting
indeterminates $\vect y=(y_0,\ldots,y_m)$ and is defined
recursively as follows. If the only edges of $G$ are loops and the
number of loops on the vertices are $e_1,\ldots,e_n$ then
\[
V_G(\vect y) = \prod_{i=1}^n y_{e_i}.
\]
Otherwise for any edge $e$ that is not a loop
\begin{equation}
V_G(\vect y) = V_{G- e}(\vect y) + V_{G/e}(\vect y).
\label{eq:Vrecursive}
\end{equation}
It is relatively simple to prove by induction that the definition
is independent of the choice of edge in~\eqref{eq:Vrecursive}.
\begin{proposition}
\[V_G(\vect y) = \EU_G\left(z_{ij} = \sum_{k=0}^j (-1)^{j-k}
\binom j k y_k\right).\]
\end{proposition}
\begin{proof}
Notice that it follows from~\eqref{eq:extW} that if the value of
$z_{ij}$ does not depend on $i$ then for any $\omega$ and
$\omega'$, $\EW_{(G,\omega)}(\vect z)= \EW_{(G,\omega')}(\vect
z)$. In particular if for all $v$, $\omega'(v)=1$, we get
$\EW_{(G,\omega)}(\vect z)= \EU_G(\vect z)$. So $\EU_G(z_{ij} =
\sum_{k=0}^j (-1)^{j-k} \binom j k y_k)$ must
satisfy~\eqref{eq:W1} and~\eqref{eq:W2}. It follows
from~\eqref{eq:W2} that $\EU_G(z_{ij} = \sum_{k=0}^j (-1)^{j-k}
\binom j k y_k)$ satisfies~\eqref{eq:Vrecursive}.
From~\eqref{eq:W1}, we see that if the only edges of $G$ are loops
and the number of loops on the vertices are $e_1,\ldots,e_n$ then
\begin{align*}
\EU_G\left(z_{ij} = \sum_{k=0}^j (-1)^{j-k}\binom j k y_k\right) &
= \prod_{i=1}^n \sum_{j=0}^{e_i} \binom{e_i}{j}
\sum_{k=0}^j (-1)^{j-k}\binom j k y_k\\
&=\prod_{i=1}^n y_{e_i}.
\end{align*}
\end{proof}

\subsection{The extended polychromate} Like the extended
$U$-polynomial, the extended polychromate is a polynomial in
countably infinitely many commuting variables $x_{i,j}$ where
$i\in \ints^{>0}$ and $j\in \ints^{\geq 0}$. It was also first
introduced in~\cite{welsh:talk05}.

The extended polychromate $\Echi$ is defined as follows.
\[
\Echi_G(\vect x) = \sum_{\pi} \Ex(\Et(\pi)),
\]
where the sum is over all partitions of $V$ and if $\Et =
((a_1,b_1),\ldots,(a_k,b_k))$ then $\Ex(\Et) = x_{a_1,b_1} \cdots
x_{a_k,b_k}$.

For example if $G$ is a triangle then
\[\Echi_G(\vect x) = x_{1,0}^3 + 3x_{2,1}x_{1,0} + x_{3,3}\]
and if $G$ is a path with two edges then
\[\Echi_G(\vect x) = x_{1,0}^3 + 2x_{2,1}x_{1,0} + x_{2,0}x_{1,0}
+ x_{3,2}.\] Note that we obtain the polychromate by substituting
$x_{i,j} = x_{i} y^{j}$ resulting in a polynomial in
$x_1,\ldots,x_n$ and $y$.

We now show that the extended polychromate and the extended Tutte
symmetric function are equivalent
\begin{theorem}\label{th:extpolysymm}
The extended polychromate and the extended Tutte symmetric
function are equivalent. More precisely for each $\Et =
((a_1,b_1), \ldots , (a_k,b_k)) \vdash (n,m)$, the coefficient of
$\Em_\Et(\vect x,\vect t)$ in $\EY_G$ is the same as the
coefficient of $\Ex_{\Et}$ in $\Echi_G$.
\end{theorem}
\begin{proof}
A colouring of $G$ induces a partition of $V$ in which two
vertices are in the same block if and only if they receive the
same colour. So we may partition the sum in~\eqref{eq:extsymmdef}
according to the partition of $V$ induced by the colouring. Hence
we can write
\[ \EY_G(\vect x;\vect t) =
\sum_{\pi}\sum_{\chi}\left(\prod_{i=1}^n x_{\chi(v_i)}\right)
\left(\prod_{i=1}^\infty (1+t_i) ^{b_i(\chi)}\right),\] where the
first summation is over all partitions of $V$ and the second over
all colourings of $V$ with strictly positive integers so that
vertices receive the same colour if and only if they are in the
same block of $\pi$. Fix a partition $\pi$ of $V$ and suppose that
$\Et(\pi) = ((a_1,b_1),\ldots,(a_k,b_k))$. Then the contribution
to $\EY$ from colourings inducing $\pi$ is $\Em_\Et(\vect x,\vect
t)$. However the monomial in $\Echi_G$ corresponding to $\pi$ is
$\Ex(\Et)$ and the result follows.
\end{proof}

\begin{corollary}
The extended polychromate and the extended $U$-polynomial are
equivalent.
\end{corollary}
\begin{proof}
This follows easily from the transitivity of equivalence.
\end{proof}
In principle one could describe a substitution in order to obtain
$\Echi$ from $\EU$ or vice versa but the procedure would be very
complicated. We show briefly how Sarmiento's result
from~\cite{sarm:poly} may be obtained as a special case of our
results by comparing the expressions linking the Tutte symmetric
function with the $U$-polynomial in Theorem~\ref{th:Usymm} and an
expression linking the Tutte symmetric function with the
polychromate deduced from Theorem~\ref{th:extpolysymm}.

Recall that the symmetric Tutte function is a function of
$(x_1,x_2,\ldots,t)$ and is a homogenous symmetric function of
degree $n$ in the $\vect x$ variables. Furthermore recall that
both the collections $\{p_\tau(\vect x) : \tau \vdash n\}$ and
$\{m_\tau(\vect x) : \tau \vdash n\}$ are bases for the homogenous
symmetric functions of degree $n$ in $\vect x$. Consequently there
are constants $a_{\tau,\tau'}$ such that $p_\tau(\vect x) =
\sum_{\tau'} a_{\tau,\tau'} m_{\tau'}(\vect x)$.

It is not difficult to compute $a_{\tau,\tau'}$. Given a partition
$\pi$, we say that the partition $\pi'$ is a coarsening of $\pi$
if every block of $\pi'$ is a union of blocks of $\pi$. Let $\pi$
be a partition of $\{1,\ldots,n\}$ of type $\tau$. Then
$a_{\tau,\tau'}$ is the number of coarsenings of $\pi$ of type
$\tau'$.

\begin{proposition}
The polychromate may be obtained from the $U$-polynomial by
replacing for each $\tau$ such that $\tau \vdash n$, the monomial
$\vect x_\tau y^j$ by $\sum_{\tau'} a_{\tau,\tau'} \vect
x_{\tau'}y^j(y-1)^{n-k(\tau)}$ where the sum is over all
$\tau'\vdash n$.
\end{proposition}
\begin{proof}
Setting $t_i=t$ for all $i$ in the extended symmetric Tutte
function we can write
\[ Y_G(\vect x,t) = \sum_{\tau \vdash n} \sum_i c_{\tau,i}m_\tau(\vect x) (1+t)^i\]
for certain constants $c_{\tau,i}$. Recall that if $\tau =
(n_1,\ldots,n_k)$ then $x_\tau = x_{n_1}\cdots x_{n_k}$. The
polychromate may be written in the form
\[ \chi_G(\vect x,t) = \sum_{\tau \vdash n} \sum_i c'_{\tau,i}{\vect x}_\tau t^i,\]
for certain constants $c'_{\tau,i}$. Theorem~\ref{th:extpolysymm}
implies that for all $\tau$ and $i$, $c_{\tau,i} = c'_{\tau,i}$.

The remarks immediately after Theorem~\ref{th:Usymm} state that
$Y_G$ may be obtained from $U_G$ by replacing the monomial $\vect
x_\tau y^i$ in $U_G$ by $p_\tau(\vect x) t^{n-k(\tau)} (t+1)^i$.
Given the relationship between the the power-sum basis and the
augmented monomial basis an equivalent substitution is to replace
$\vect x_\tau y^i$ by $t^{n-k(\tau)}(t+1)^i \sum_{\tau'\vdash n}
a_{\tau,\tau'} m_{\tau'}(\vect x)$.

Now the first part of the proof shows that replacing $t$ by $y-1$
and $m_\tau(\vect x)$ by ${\vect x}_\tau$ in $Y_G$ gives $\chi_G$
and the result follows.
\end{proof}
A similar argument shows how to obtain $U$ from the polychromate.

\section{Conclusions and open problems}

The graph polynomials and symmetric functions that we have
discussed are related by the following partial order where a
function $P$ is above $Q$ if $P$ specializes to $Q$.

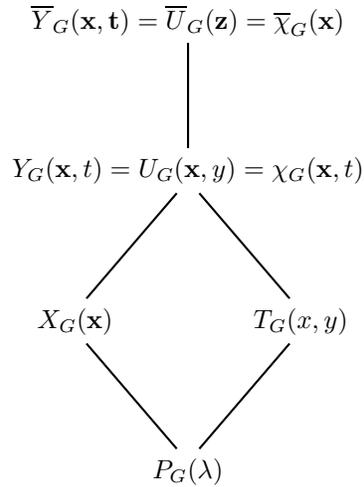
\begin{center}
\pspicture*(-5,-0.5)(5,6.5) \rput(0,0){$P_G(\lambda)$}
\rput(-1.5,2){$X_G(\vect x)$}
\rput(1.5,2){$T_G(x,y)$}

\rput(0,4){$Y_G(\vect x,t)=U_G(\vect x,y)=\chi_G(\vect x,t)$}
\rput(0,6){$\EY_G(\vect x,\vect t)=\EU_G(\vect z)=\Echi_G(\vect x
)$}

\psline(0,5.7)(0,4.3)
\psline(1.35,2.3)(0.15,3.7)\psline(-1.35,2.3)(-0.15,3.7)
\psline(1.35,1.7)(0.15,0.3)
\psline(-1.35,1.7)(-0.15,0.3)

\endpspicture
\end{center}
The relationships between many other polynomials are considered
in~\cite{mak:zoo}.

An open problem is to find a pair of loopless graphs $G_1$, $G_2$
for which $U_{G_1}=U_{G_2}$ (or for which either of the other
equivalent functions coincide) but $\EU_{G_1}\ne\EU_{G_2}$. The
following graphs are the smallest known pair of non-isomorphic
graphs with the same polychromate~\cite{bry:intersection}. However
it is easy to see that they also have the same extended
$U$-polynomial.

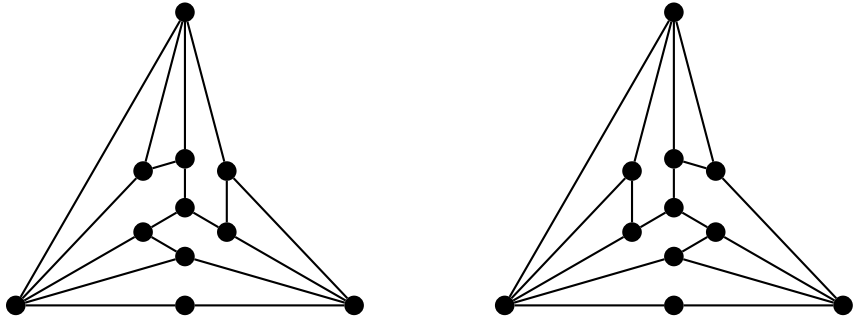
\begin{center}
\psset{unit=1.3} \pspicture(-2,-2)(7,2.2)

\cnode*(0,0){0.1}{a} \cnode*(0,2){0.1}{b}
\cnode*(-1.73,-1){0.1}{c} \cnode*(0,-1){0.1}{d}
\cnode*(1.73,-1){0.1}{e} \cnode*(0,-0.5){0.1}{f}
\cnode*(0,0.5){0.1}{g} \cnode*(-0.428,-0.25){0.1}{h}
\cnode*(0.428,-0.25){0.1}{i} \cnode*(-0.428,0.375){0.1}{j}
\cnode*(0.428,0.375){0.1}{k}

\ncLine{-}{a}{g}\ncLine{-}{a}{h}\ncLine{-}{a}{i}
\ncLine{-}{b}{c}\ncLine{-}{b}{j}\ncLine{-}{b}{k}
\ncLine{-}{c}{d}\ncLine{-}{f}{c}\ncLine{-}{d}{e}
\ncLine{-}{g}{j}\ncLine{-}{h}{f}\ncLine{-}{i}{k}
\ncLine{-}{b}{g}\ncLine{-}{c}{j}\ncLine{-}{e}{k}
\ncLine{-}{e}{f}\ncLine{-}{c}{h}\ncLine{-}{e}{i}

\cnode*(5,0){0.1}{aa} \cnode*(5,2){0.1}{bb}
\cnode*(3.27,-1){0.1}{cc} \cnode*(5,-1){0.1}{dd}
\cnode*(6.73,-1){0.1}{ee} \cnode*(5,-0.5){0.1}{ff}
\cnode*(5,0.5){0.1}{gg} \cnode*(4.572,-0.25){0.1}{hh}
\cnode*(5.428,-0.25){0.1}{ii} \cnode*(4.572,0.375){0.1}{jj}
\cnode*(5.428,0.375){0.1}{kk}

\ncLine{-}{aa}{gg}\ncLine{-}{aa}{hh}\ncLine{-}{aa}{ii}
\ncLine{-}{bb}{cc}\ncLine{-}{bb}{jj}\ncLine{-}{bb}{kk}
\ncLine{-}{cc}{dd}\ncLine{-}{ff}{cc}\ncLine{-}{dd}{ee}
\ncLine{-}{hh}{jj}\ncLine{-}{ii}{ff}\ncLine{-}{gg}{kk}
\ncLine{-}{bb}{gg}\ncLine{-}{cc}{jj}\ncLine{-}{ee}{kk}
\ncLine{-}{ee}{ff}\ncLine{-}{cc}{hh}\ncLine{-}{ee}{ii}

\endpspicture
\end{center}
\section*{Acknowledgements}
We would like to thank Dominic Welsh and Janos Makowsky for useful discussions.


\end{document}